\documentclass{article}
\usepackage[utf8]{inputenc}
\usepackage{authblk}
\providecommand{\keywords}[1]{\textbf{\textit{Keywords---}} #1}

\usepackage[a4paper, total={6in, 8in}]{geometry}
\usepackage{chemformula}

\usepackage{graphicx}
\usepackage[breaklinks=true,colorlinks=true,linkcolor=blue,urlcolor=blue,citecolor=blue]{hyperref}

\usepackage{amsfonts}
\usepackage{amsmath}
\usepackage{multirow}
\usepackage{tikz}
\tikzstyle{mybox} = [draw=black, thick,
 rectangle, rounded corners, inner sep=10pt, inner ysep=20pt]
\tikzstyle{fancytitle} =[fill=black, text=white]

\renewcommand{\bf}[1]{\boldsymbol{#1}}
\renewcommand{\d}{\mathrm{d}}
\newcommand{\bx}{\boldsymbol{x}}
\newcommand{\by}{\boldsymbol{y}}

\usepackage{booktabs,siunitx} 
\usepackage{multirow}

\usepackage[boxed,ruled,vlined]{algorithm2e}

\SetCommentSty{mycommfont}

\title{Revealing hidden dynamics from time-series data by ODENet}
\author[ ]{Pipi Hu\textsuperscript{1$*$},Wuyue Yang$^{1}$\thanks{These authors have contributed equally to this work}}
\author[ ]{Yi Zhu$^{1}$\thanks{yizhu@mail.tsinghua.edu.cn}}
\author[ ]{Liu Hong$^{2}$\thanks{hongliu@sysu.edu.cn}}

\affil[1]{Yau Mathematical Sciences Center, Tsinghua University, Beijing, 100084, China.}

\affil[2]{School of Mathematics, Sun Yat-Sen University, Guangzhou, 510275, China.}

\begin{document}
\maketitle

\begin{abstract}
To derive the hidden dynamics from observed data is one of the fundamental but also challenging problems in many different fields. In this study, we propose a new type of interpretable network called the ordinary differential equation network (ODENet), in which the numerical integration of explicit ordinary differential equations (ODEs) are embedded into the machine learning scheme to build a general framework for revealing the hidden dynamics buried in massive time-series data efficiently and reliably. ODENet takes full advantage of both machine learning algorithms and ODE modeling. On one hand, the embedding of ODEs makes the framework more interpretable benefiting from the mature theories of ODEs. On the other hand, the schemes of machine learning enable data handling, paralleling, and optimization to be easily and efficiently implemented.  From classical Lotka-Volterra equations to chaotic Lorenz equations, the ODENet exhibits its remarkable capability in handling time-series data even in the presence of large noise. We further apply the ODENet to real actin aggregation data, which shows an impressive performance as well. These results demonstrate the superiority of ODENet in dealing with noisy data, data with either non-equal spacing or large sampling time steps over other traditional machine learning algorithms.
\end{abstract}

\keywords{ODENet, time-series data, ordinary differential equations, chemical reactions}

\section{Introduction}
At every moment, massive data has been collected through diverse human activities. 
 And revealing the hidden dynamics from those collected time-series data is one fundamental goal of science.
There are many ``standard" theories to describe such dynamics, among which differential equations are probably the most successful one. For example, Newton's equation $F = m\ddot{x}$ combined with the law of universal gravitation gives a simple and correct picture to explain the complicated motions of planets and the sun and even the existence of Pluto more than half a century before its discovery.  However, many new fields, like psychology and social science, still lack rigorous and quantitative theories/equations until now. Therefore, how to construct effective models from a data-driven point of view becomes an interesting topic.

Unfortunately, the time-series data in record usually contain a lot of missing points and even flaws. They are also highly noisy, with useful signals deeply buried. These facts make analyzing time-series data and extracting useful models or principles hard and tricky. A most famous example is the explanation of planetary orbits in the solar system. Even though wrongly placing the earth at the center, Claudius Ptolemy was still able to explain the motion of planets and the sun by obscure deferent and epicycle with certain accuracy. Only 1,400 years later, after the landmark works of Copernicus, Kepler, and Newton, a much simpler and more correct picture could be gradually established and accepted. This story highlights the ambiguity and difficulty behind data-driven modeling. 

In recent years, the analysis of time-series data has already become a specific subject \cite{little2013oxford}. Especially in the presence of big data, plenty of machine learning algorithms, like recurrent neural network (RNN) \cite{elman1990finding}, long short-term memory (LSTM) \cite{hochreiter1997long}, \textit{etc.}, have been widely applied to various fields. RNN, which uses its internal state network to store representations of recent inputs and reuse the output to process the time series, shows a dramatic different ``memory'' property and network flow structure from other supervised neural networks.  LSTM solves the problems of vanishing gradients in RNN by using the so-called ``long term'' and ``short term'' memories. Though RNN and LSTM are very successful in applications, their performance in dealing with time-series data collected from physical processes are not very promising \cite{chen2018neural}. By adding short-cuts to jump over some layers, the ResNet\cite{he2016deep} can avoid the problem of vanishing gradients and shows better performances than classical deep learning networks. Those short-cuts can also be treated as some kind of ``memory'', which keeps the intrinsic properties unchanged during the learning procedure.  Mathematically, the ResNet is analogous to the numerical schemes of ODEs. This interesting connection leads to the so-called continuous view of machine learning, which has been explored a lot from a mathematical point of view \cite{weinan2017proposal,ma2019machine} and also been used for constructing new machine learning algorithms \cite{raissi2018multistep, chen2018neural}.

Besides the above mentioned neural networks, there are many other attempts to extract dynamical equations from the time-series data in the past years. For example, Bongard and Schmidt used symbolic regression to find nonlinear differential equations \cite{bongard2007automated,schmidt2009distilling,koza1992genetic}. Kutz \textit{et al.} proposed a framework named ``SINDy'' by combining regression and sparse identification to reveal nonlinear dynamical systems \cite{kutz16sparse}. To stabilize the performance of sparse identification in the presence of noise, some technical schemes for differential operations were proposed \cite{rudy2017data, kutz19noise}. Recently, Dong \textit{et al.} \cite{long2018pde,long2019pde} use kernels in CNN which mimic the differential operators to reveal the PDE dynamics from the training data.

Inspired by the great success of modern machine learning algorithms, in the current study we aim to reveal the hidden dynamics from the time-series data without prior knowledge. Different from most previous works that focus on the efficiency and accuracy of predictions of neural networks with little or no physical understandings, we are more interested in deriving the explicit governing equations from the time-series data, which is done under the help of a new type of network, called the ordinary differential equation network (ODENet). By combining  the optimization structure of neural networks with symbolic regression, sparse identification, and signal-noise decomposition,
ODENet exhibits an outstanding ability in deriving the explicit ODE models for population dynamics modeled by Lotka-Volterra equations in presence of large noise, strange-attractors of Lorenz equations in the chaotic region, as well as the hidden actin growth dynamics and molecular mechanisms base on real experimental data under distinct conditions. In particular, through these studies our framework has been proven to be robust, noise-tolerant, immune to unequal time steps of training data, and therefore ODENet is quite suitable for time-series data analysis and data-driven mathematical modeling.

\section{The architecture of ODENet}

There are two ways to interpret the dynamics behind the time-series data. The usual machine learning algorithms, like LSTM and deep learning, tend to use a vast neural network containing a large number of free parameters to achieve the goal of representing a complex mapping function that best fits the data set. Through iterative training and optimization, the data correlation is transformed into very complicated and thus unexplainable relations among network nodes. In contrast, regression and sparse identification methods adopt an alternative view, which is more concerned about the construction of explicit relations or dynamic differential equations for a globally fitting of the data but only with a few parameters. Each way has its advantages and appropriate applicable regions. Since we are more interested in the physical mechanisms and mathematical models behind the data, the parts of the two ways are combined and realized through the ODENet in the current study.

\begin{figure}[!htp]
\centering
\includegraphics[width=\linewidth]{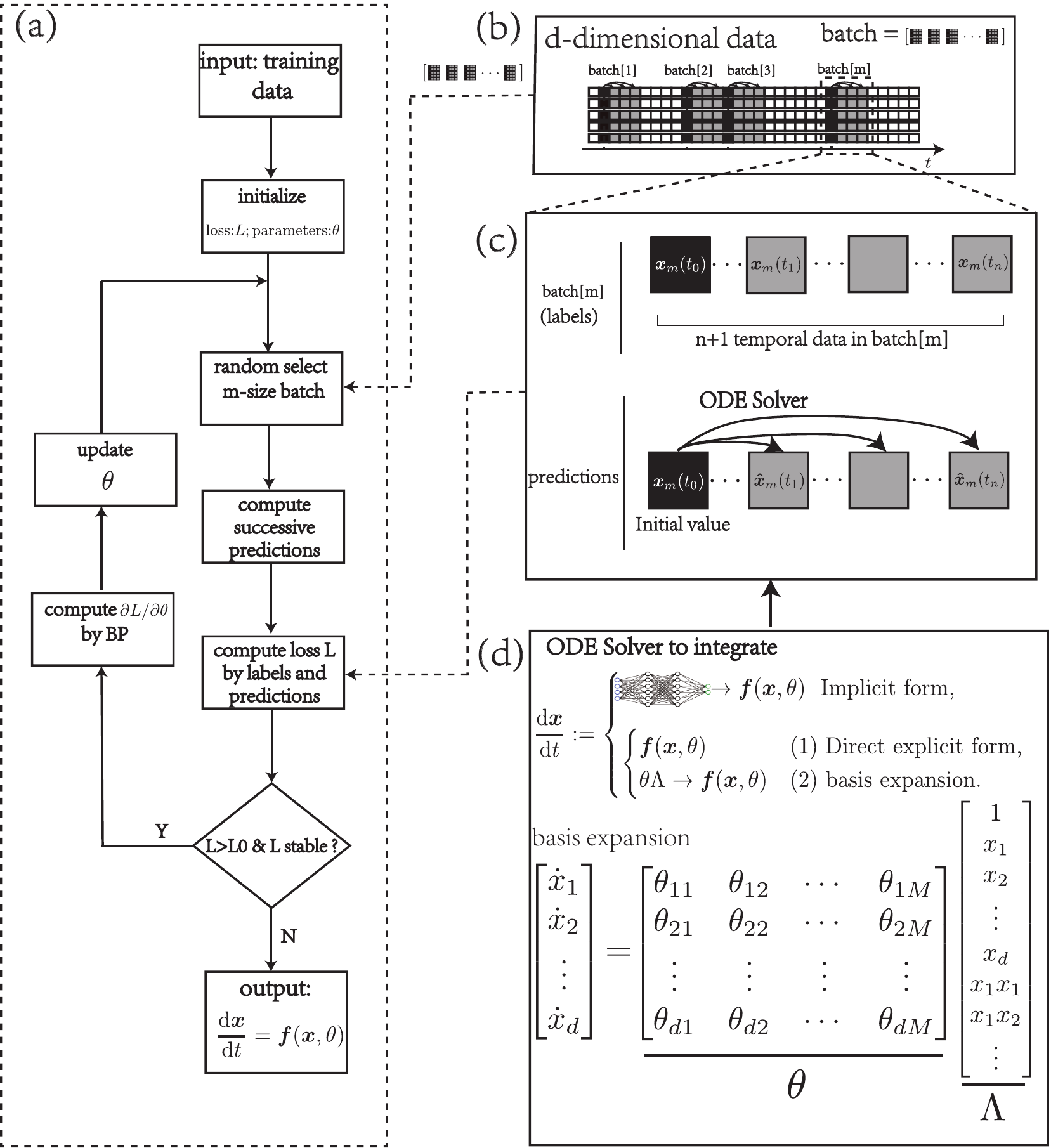}
\caption{The architecture of ODENet with key steps highlighted in the right columns.  (a) shows the basic flowchart of the training algorithm, (b) exhibits the structure of training data, and (c) rooms in to show the data structure of one training piece in one batch, and (d) gives the mathematical schemes of ODEs.}
\label{fig.ODENet}
\end{figure}

As a modification of ResNet, the basic structure of ODENet (see Figure \ref{fig.ODENet}) mimics the numerical solvation of ODEs with unspecified parameters to be learned from the data. Through iteratively minimizing the difference between predicted time trajectories and training data measured by a certain loss function, unknown parameters are optimized in such a way that the most suitable ODE model for the given time-series data is explicitly specified.

 \textit{\textbf{Data Batching:}}   To be concrete,
the learning procedure of the ODENet begins with $m$ randomly selected points from the training data set, $\bf{x}_i(t_0)$, for $i=1,2,\cdots, m$, which will be used as the starting points for  integration. For each starting point, its following $n$ successive data are picked as labels. Till now, we have extracted $m$ pieces of time-series data of length $n+1$, \textit{i.e.} $\bf{x}_1(t_0), \bf{x}_1(t_1), \cdots \bf{x}_1(t_n)$; $\bf{x}_2(t_0), \bf{x}_2(t_1), \cdots \bf{x}_2(t_n)$; $\cdots$, $\bf{x}_m(t_0), \bf{x}_m(t_1), \cdots \bf{x}_m(t_n)$, which constitute the batch shown in Figure \ref{fig.ODENet}b. In general, $\bf{x}_i(t_j)$ is a $d$-dimensional vector, where index $i$ represents the $i^{th}$ piece and index $j$ represents the $(j+1)^{th}$ point in the time-series data counting from the starting one. Note $t_i-t_{i-1}$ may not necessarily be the same for different $i\in\{0,1,2,..,n\}$. However, as large $\Delta t$ will make the problem stiff while small $\Delta t$ may be unable to provide sufficient information for the dynamics, the time steps must be carefully adjusted in order to keep a well balance between model accuracy and stiffness.

\textit{\textbf{ODE Dynamics:}}  Next, referring to each piece of data in the batch, we suppose they satisfy the following initial-value problem of a system of autonomous ordinary differential equations,
\begin{equation}\label{eq.ode}
\begin{split}
&\frac{\d \bx}{\d t} = \bf{f}(\bx;\bf{\theta}),\\
&\bx{(t_0)} = (x_1{(t_0)},x_2{(t_0)},\cdots,x_d{(t_0)})^T.\\
\end{split}
\end{equation}
This is the key assumption of our ODENet. Here, the starting point $\bx{(t_0)}$ at $t_0$ is taken as the initial value. $\bf{\theta}$ stands for the unspecified parameters, which determine the explicit form of the right-hand side terms, and will be addressed later. Next, above ODEs will be solved numerically by mutual ODE solvers, like the Runge-Kutta method, whose numerical solutions at $t_1,\cdots, t_n$ are denoted as $\bf{\hat{x}}(t_1), \bf{\hat{x}}(t_2), \cdots \bf{\hat{x}}(t_n)$ respectively, in order to distinguish from the labels $\bf{x}(t_1), \bf{x}(t_2), \cdots \bf{x}(t_n)$. 

Therefore, we actually adopt an integral method to try to predict the ODE dynamics starting from $\bx{(t_0)}$ at $t_0$. Unlike most previous regression based methods, \textit{e.g.} SINDy \cite{kutz16sparse}, our approach minimizes the deviations from the training data $\bf{x}(t_1), \bf{x}(t_2), \cdots \bf{x}(t_n)$ at time $t_1, t_2, \cdots, t_n$ instead of the  derivatives. We will come back to this point soon.

\textit{\textbf{Approximation of Unknown Functions:}} 
Before turning to the description of the loss function and optimization scheme, we emphasis here the structure of the right-hand side terms $\bf{f}(\bx;\bf{\theta})=(f_1,f_2,\cdots,f_d)^T$, whose form completely determines the ODE dynamics in \eqref{eq.ode} and has to be found out from the training data.

There are two different ways -- the implicit and the explicit. The usual neural networks, especially deep neural networks, which show a great ability to approximate very complicated mapping functions, in reality, belong to the former since they often contain multiple layers, tremendous nodes, and connections with massive adjustable free parameters. In most cases, they are ``black boxes'' to us. The implicit approach has been adopted by several groups previously and applied to MINIST as an example \cite{chen2018neural,raissi2018multistep}.

In the current study, as we are more interested in the dynamics and mechanism behind the time-series data, we hope to derive the ODE model in an explicit way.  To achieve such a goal, the right-hand side terms $\bf{f}(\bx;\bf{\theta})$ are expanded through basis functions with corresponding coefficients $\bf{\theta}$. Although the selection of proper basis function is quite tricky and problem dependent, polynomials are among the most often used ones in practice. For the $d$-dimensional vector $\bx = (x_1, x_2, \cdots, x_d)^T$, the $p^{th}$-order complete polynomials 
 $$\{1, x_1, x_2, \cdots, x_d, x_1x_1, x_1x_2,\cdots, x_d^p\}$$ 
  have $M={p+d \choose p}$ terms. So that we have
\begin{equation}
\bf{f}(\bx;\bf{\theta}) = \bf{\theta\Lambda},
\label{eq.right}
\end{equation}
where $\bf{\Lambda} = (1, x_1, x_2, \cdots, x_d, x_1x_1, x_1x_2,\cdots, x_d^p)^T$ is the complete set of $p^{th}$-order polynomial basis with coefficients $\bf{\theta}=(\theta_{ij})_{d\times M}$.

\textit{\textbf {Loss Function:}} By tuning free parameters $\bf{\theta}$, which specify the concrete form of ODEs, the loss function
\begin{equation}
\mathcal{L(\bf{\theta};\lambda)} = \|\bx-\tilde{\bx}\|_2+\mu\|\bf{\theta}\|_1
\end{equation}
is expected to be minimized.  Here the first term characterizes the difference between the training data and predictions, while the second term represents the sparsity requirement with $\mu>0$ as a hyper parameter. According to the Occam's Razor, ``plurality should not be posited without necessity.''  So it is expected that the components of $\bf{\theta}$ should be as many zeros as possible, which thus corresponds to the simplest model and also the smallest $L1$ norm.  Here we may encounter the
$L1$ optimization problem. The autograd adopted in PyTorch uses subderivatives for calculating gradients of $L1$ norm at the non-differentiable points. We refer readers to alternative new algorithms, such as the Split Bregman Method \cite{goldstein2009split}, Coordinate Descent \cite{wright2015coordinate}, Proximal Gradient Method \cite{schmidt2011convergence} and \textit{etc.},
which may offer a better performance on optimizing $L1$ norm than the subgradients.

Besides $\mu$, and additional threshold parameter $\gamma>0$  is adopted to accelerate the learning of sparse models. Once a component of $\theta$ is smaller than $\gamma$, it will be forced to be zero in the remain learning process. As a consequence, $\gamma$ should be carefully tuned to make a good balance between model simplicity and fitting accuracy. Generally speaking, within a certain region, the larger the threshold is, the faster the loss function will converge. However, once the threshold becomes too large, the learning process has a very high risk of failing due to too many abandoned terms.

\textit{\textbf {Parameter Optimization:}} The searching of the best $\bf{\theta}$ is a global optimization problem, which can be done via mature learning algorithms in the neural network, for example, gradient descent or stochastic gradient descent combined with backward propagation algorithms \cite{lecun1989backpropagation}. Especially in PyTorch, this can be simply done by an integrated autograd library. Once the forward flow is constructed, the backward flow will be automatically built by Pytorch \cite{paszke2017automatic}. It should be noted that autograd is an automatic differentiation method, whose computational graph will become too large and make the problem computationally very expensive once we want to simulate the system for a long time. To avoid this difficulty, we have to look into the problem of stiffness, which is quite often encountered in numerical simulations. Finally, repeating the above procedure iteratively until the loss function does not decay efficiently anymore or less than a threshold, we finish the learning procedure of ODENet, which is summarized through pseudocodes as below. \\

\begin{algorithm}[H]
\SetAlgoLined
\KwIn{time-series data $\bf{x}(0),\bf{x}(t_1),\cdots, \bf{x}(t_n)$, where $\bf{x}=(x_1,x_2,\cdots, x_d)^T$.}
\KwOut{d-dimensional first order ODEs $\frac{\d\bf{x}}{\d t}=\bf{f}(\bf{x},\theta)$.}
initialize parameters $\theta$\;
initialize hyperparameters: $threshold\_L$, $threshold\_\theta$, $m$,$n$\;

 \While{$L>threshold\_L$}{
 \tcp*[h]{Construct one batch}\\
 select $n$-length intervals from $m$ random positions for a batch\;
$batch\leftarrow [[\bx_0^1,\bx_1^1,\cdots,\bx_n^1],[\bx_0^2,\bx_1^2,\cdots,\bx_n^2],\cdots,[\bx_0^m,\bx_1^m,\cdots,\bx_n^m]]$\;
 $batch\_t\leftarrow [[t_0^1,t_1^1,\cdots,t_n^1],[t_0^2,t_1^2,\cdots,t_n^2],\cdots,[t_0^m,t_1^m,\cdots,t_n^m]]$\;
 $batch\_init\leftarrow batch[:,1]$\;
 $batch\_label\leftarrow batch[:,2:]$\;
 $batch\_size\leftarrow m$\;
 $L\leftarrow 0$\;
 \For{$i\leftarrow $ 1 \KwTo $batch\_size$}{
 \tcp*[h]{Compute the predictions based on $batch\_init$}\
$pred[i,:]\leftarrow ODESolve(batch\_init[i],\,batch\_t[i,:],\,\theta)$\;
 \tcp*[h]{The $\theta$ matrix should be sparse matrix}\\
$L\leftarrow L + \|Abs(pred[i,:]-batch\_label[i,:])\|_2+\mu \|\theta\|_1$\;
}
$\theta_{\cdot}grad\leftarrow \frac{\partial L}{\partial\theta}$ by BP algorithm\;
Update parameters $\theta$ by Adam methods\cite{kingma2014adam};
 \tcp*[h]{ Set $\theta_{ij}$ which is small enough to be zero}\\
$\theta[\theta<\gamma]\leftarrow 0$\;
}
 \caption{Pesudocode of ODENet (PyTorch)}
\end{algorithm}

\textit{\textbf{Stiffness Problem:}}  It is noted that we put no constraint during the step of data selection. On one hand, this provides great facility in dealing with real data; on the other hand, it also leaves us in danger of facing stiffness problems from time to time during numerical simulations, which becomes especially outstanding when linear combinations of high-order polynomials are taken as the right-hand side terms of an ODE system. 

According to the general knowledge of ODE numerics, we suggest the following ways for dealing with the stiffness problem.  First, divide long time trajectories of the training data set into many short pieces. And thus each piece will contain only a few time steps, which could be easily solved by classical gradient methods implemented in PyTorch. Second, towards the ODE solver, we suggest using self-adaptive and/or implicit ODE solvers, such as dopri5 which was designed for solving stiff cases for the best \cite{calvo1990fifth,hairer1999stiff}. Third, instead of arbitrarily setting parameters (coefficients of the polynomial basis) through random number generators, regression methods can be introduced to make a reasonable initial guess on the parameters  \cite{kutz16sparse}. Last, for recent advances in the direction of multiscale modeling methods combined with ODE or PDE based neural networks, which may provide an alternative solution to the stiffness problems, see \textit{e.g.} Ref. \cite{han2019uniformly, Yang2020When}.

\textit{\textbf {Integral v.s. Differential:}} 
One significant feature of the ODENet from previous regression based methods like SINDy \cite{kutz16sparse} is that our approach is based on an integration of explicit ODEs along the  time trajectory, while theirs \cite{bongard2007automated,schmidt2009distilling,kutz16sparse} are all based on differentiation between neighboring points. To be concrete, in our approach the data points on a time trajectory are predicted based on an integral solution of the ODE model, i.e. $\bx(t)=\int_{t_0}^{t_n}\bf{f}(\bx(\tau);\bf{\theta})d\tau+\bx(t_0)$. The loss function is designed to minimize the difference between the predicted $\bx(t)$ and real ones. In contrast, according to the methods reported in Refs. \cite{bongard2007automated,schmidt2009distilling,kutz16sparse}, what they actually tried to minimize is the difference between the function $\bf{f}(\bx(\tau);\bf{\theta})$ on the right-hand side of \eqref{eq.ode} and the time derivative 
$\dot{\bx}(t)=\frac{\bx(t)-\bx(\tau)}{t-\tau}$ calculated from the training data. Clearly, the latter is a kind of differentiation methods.

The differentiation methods are simple and straightforward, easy for implementation, numerically fast and efficient, and suitable for high-dimensional complicated dynamics. In comparison, the integral methods enjoy advantages like more stable against faults and flaws in data, more noise-tolerant, able to endure data with large time steps, \textit{etc.}

\textit{\textbf {Noise:}} As we have discussed in the introduction, real data may have noise, flaws and faults. To deal with this issue, we take a brutal method here by incorporating the strength of noise as learning parameters too. Suppose there is a finite time series $\by(t)\in \mathbb{R}^d$ with noise $\bf{e}(t) = \epsilon\|\by\|_{\infty}\bf{\eta}$, where $\epsilon$ denotes the noise strength, $\|\by\|_{\infty}$ is the maximal value in $y$, and $\bf{\eta} \sim \bf{\mathcal{N}}(0,1)$ are d-dimensional normally distributed random variables. ODENet is proposed to extract the following autonomous system of ODEs
\begin{equation}
\frac{\d \bx}{\d t} = \bf{f}(\bx;\bf{\theta}), ~ \bx(t) = \by(t)-\bf{e}(t)\in\mathbb{R}^d.
\label{eq.model}
\end{equation}
And we can apply the same process as before to extract the governing dynamics from the data. Apparently, at this time, the loss function depends on $\bf{e}(t)$ too, i.e. $L=L(\bf{\theta},\bf{e}(t);\mu)= \|\hat{\bx}+\hat{\bf{e}}-\by\|_2+\mu\|\bf{\theta}\|_1$, where $\hat{\bf{x}}$ and $\hat{\bf{e}}$ denote the output of ODENet and learned noise respectively.  Although directly incorporating noise as learning parameters may cause a big increment in the computational cost, it allows the treatment of large noise or color noise in principle. If the noise is relatively small, it is more appropriate to try to learn the original data directly and then check its robustness against noise. Interested readers may refer to Refs. \cite{rudy2017data, kutz19noise} for alternative solutions.

\section{Numerical experiments}\label{sec.num}

In this section, we are going to apply the ODENet to the study of Lotka-Volterra equations in diverse parameter regimes and Lorenz equations in the chaotic regime. Through these examples, the power and advantage of ODENet could be demonstrated. 

\subsection{Lotka-Volterra equations with and without large noise}

The Lotka-Volterra (LV) equations, also known as predator-prey equations, were first introduced by Lotka  \cite{lotka1926elements} and Volterra \cite{volterra1927variazioni} in the 1920s to describe the population dynamics of preys interacting with predators in ecological systems. LV equations have been widely applied to ecological balance \cite{samuelson1971generalized}, environmental protection \cite{tsai2016elucidating}, disease prevention and control \cite{holt1985infectious}, etc. A very general form of LV equations including species growth and death, intraspecies and interspecies competition reads
\begin{equation}\label{eq.lv}
\begin{split}
\frac{\d x_1}{\d t} = C_{11} x_1 + C_{12} x_1x_2 + C_{13} x_1^2,\\
\frac{\d x_2}{\d t} = C_{21} x_2 + C_{22} x_1x_2 + C_{23} x_2^2.\\
\end{split}
\end{equation}

It's easy to see above equations have four fixed points, i.e.
$\left( {0, 0} \right)$, $\left( { - \frac{{{C_{11}}}}{{{C_{13}}}},0} \right)$, $\left( {0, - \frac{{{C_{21}}}}{{{C_{23}}}}} \right)$ and $\left( { - \frac{{{C_{12}} * {C_{21}} - {C_{11}} * {C_{23}}}}{{{C_{12}} * {C_{22}} - {C_{13}} * {C_{23}}}}, - \frac{{{C_{11}} * {C_{22}} - {C_{13}} * {C_{21}}}}{{{C_{12}} * {C_{22}} - {C_{13}} * {C_{23}}}}} \right)$ under the condition ${C_{12}} \times {C_{22}} \ne {C_{13}} \times {C_{23}}$,${C_{13}} \ne 0$ and ${C_{23}} \ne 0$. And the dynamic behaviors of LV equations around these fixed points $(x_1^*,x_2^*)$ are fully specified by the Jacobian matrix (its eigenvalues to be exact)
\begin{equation}\label{eq.jaco}
\begin{split}
J = {\left[ {\begin{array}{*{20}{c}}
 {{C_{11}} + {C_{12}}{x_2} + 2{C_{13}}{x_1}}&{{C_{12}}{x_1}}\\
 {{C_{22}}{x_2}}&{{C_{21}} + {C_{22}}{x_1} + 2{C_{23}}{x_2}}
 \end{array}} \right]_{\left( {x_1^*,x_2^*} \right)}},
\end{split}
\end{equation}
which can be roughly classified into three basic types -- the extinction of one species (over damped), or the evolution to an equilibrated coexistence (spiral), or to a continuing oscillation (limit cycle) \cite{may1972limit}.

\begin{figure}[!htp]
\centering
\includegraphics[width=\linewidth]{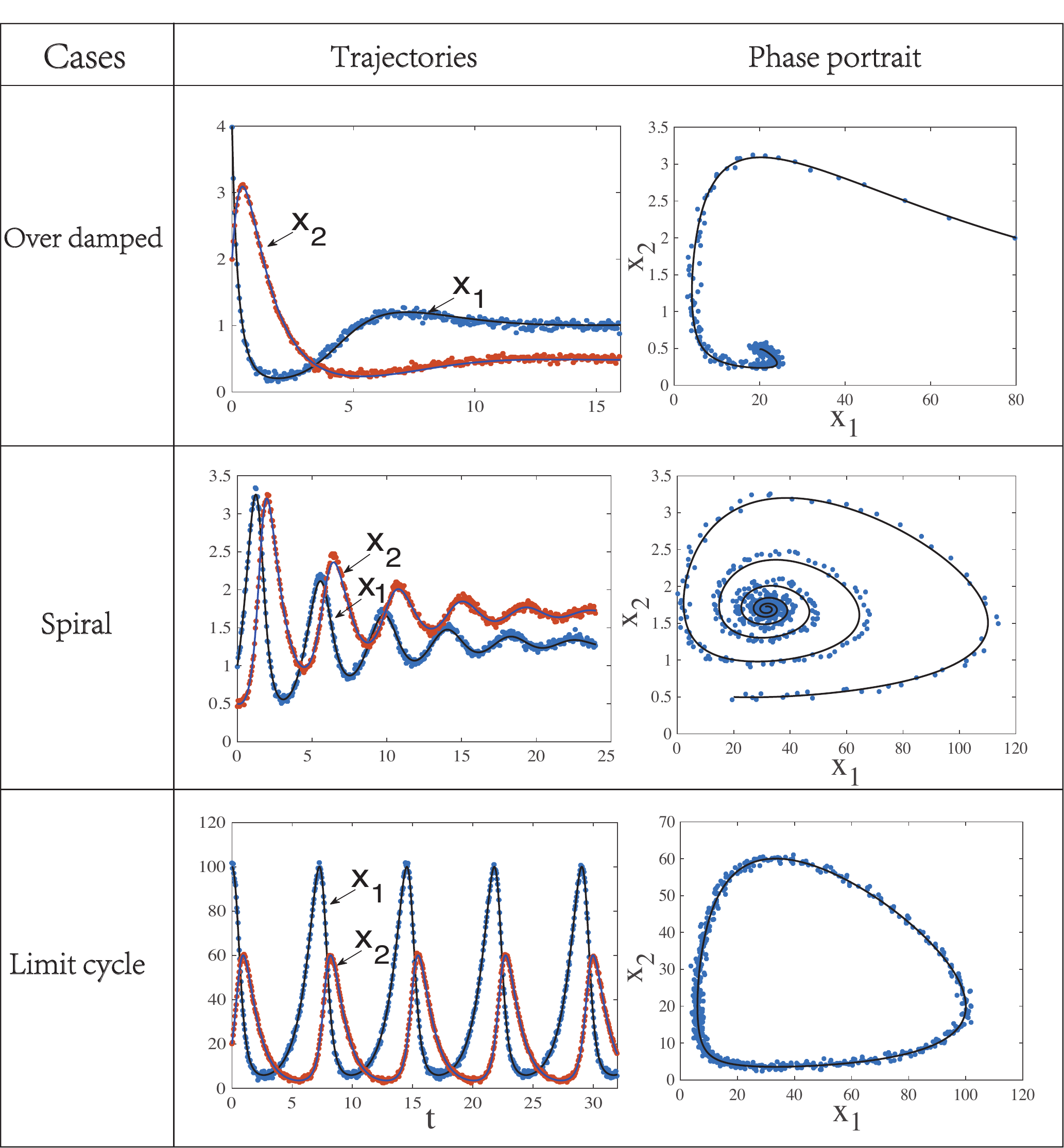}
\caption{The accuracy of ODENet predictions in both time domain (left column) and phase space (right column) in comparison with exact solutions of LV equations.  1\% white noise is added.}
\label{fig.lv}
\end{figure}

\begin{table}[htp]
\small
\setlength{\tabcolsep}{5pt}
\centering
\renewcommand{\arraystretch}{1.2}
\begin{tabular}{@{} l l cccccc @{}} 
\toprule
\toprule

\multicolumn{2}{c}{\multirow{2}*{LV equations}} & \multicolumn{6}{c@{}}{ \hspace{0em} Parameters}\\ 
\cmidrule(l){3-8}
  & & $C_{11}$& $C_{12}$ & $C_{13}$& $C_{21}$& $C_{22}$& $C_{23}$ \\
  
\midrule
\multirow{4}{4.5em}{ Over damped 1\% noise}& Model & 1.5& -1& -1& -1& 1& 0\\\cmidrule(l){2-8}
& ODENet &1.49& -9.87$\times 10^{-1}$&-9.94$\times 10^{-1}$& -1.00& 9.95$\times 10^{-1}$& 0\\
&SINDy$\sharp1$ &---&---&---&---&---&---\\ 
&SINDy$\sharp2$ &1.50&-1.00&-1.00&-1.00  &1.00&0\\\hline 
\multirow{4}{4.5em}{ Spiral 1\% noise}& Model& 2& -1.1& -0.1& -1& -0.1& 0.9 \\\cmidrule(l){2-8}
&ODENet &1.99& -1.11& -9.91$\times 10^{-2}$& -9.87$\times 10^{-1}$& -1.06$\times 10^{-1}$& 9.02$\times 10^{-1}$ \\
&SINDy$\sharp1$ &---&---&---&--- &---&---\\ 
&SINDy$\sharp2$ &2.18&-8.54$\times 10^{-1}$&-1.04$\times 10^{-1}$&-9.89$\times 10^{-1}$&-1.01$\times 10^{-1}$&8.65$\times 10^{-1}$\\\hline
\multirow{4}{4.5em}{Limit cycle \\1\% noise}& Model &1& -0.05& 0& -1& 0.03& 0 \\\cmidrule(l){2-8}
 & ODENet& 9.97$\times 10^{-1}$ &-4.98$\times 10^{-2}$ &0 &-1.00 &2.99$\times 10^{-2}$& 0\\
 &SINDy$\sharp1$ &1.00&-4.98$\times 10^{-2}$&0&-9.98$\times 10^{-1}$  &3.01$\times 10^{-2}$&0\\ 
&SINDy$\sharp2$ &1.01&-4.99$\times 10^{-2}$&0&-9.95$\times 10^{-1}$&3.00$\times 10^{-2}$&0\\\hline
\multirow{4}{4.5em}{Limit cycle\\10\% noise}& Model &1& -0.05& 0& -1& 0.03& 0 \\\cmidrule(l){2-8}
& ODENet& 9.69$\times 10^{-1}$ &-4.89$\times 10^{-2}$ &0 &-9.82$\times 10^{-1}$ &2.96$\times 10^{-2}$& 0\\
 &SINDy$\sharp1$ &9.83$\times 10^{-1}$&-5.04$\times 10^{-2}$&0&-1.00  &3.02$\times 10^{-2}$&0\\ 
&SINDy$\sharp2$&1.00&-4.99$\times 10^{-2}$&0&-9.97$\times 10^{-1}$&3.00$\times 10^{-2}$&0\\

\bottomrule
\bottomrule
\end{tabular}
\caption{Comparison of ODENet and SINDy on model coefficients for three typical dynamics of LV equations. For comparison, different data sampling time steps are adopted. For ODENet and SINDy$\sharp1$, $\Delta t = 0.01$, and for SINDy$\sharp2$, $\Delta t = 0.001$. Instead of 12 coefficients for a 2-order complete polynomial basis, only $6$ coefficients corresponding to those in \eqref{eq.lv} are listed for simplicity. No redundant coefficients have been learned by ODENet, in contrast to SINDy. The symbol ``---" indicates the failing in deriving a reasonable model from the given data set.   In this study, SINDy has been performed under different hyper-parameters for 50 independent runs. The best result is picked out and shown in the table.}
 \label{tab.lv}
\end{table}

\begin{table}[!htp]
\scriptsize
\setlength\extrarowheight{2pt}
\setlength{\tabcolsep}{13.5pt}
\centering
\begin{tabular}{@{} l ccccccc @{}} 
\toprule
\toprule
\multirow{2}*{Limit cycle} & \multicolumn{6}{c@{}}{ \hspace{6em} Parameters}\\
\cmidrule(l){3-8}
 & &$C_{11}$& $C_{12}$ & $C_{13}$& $C_{21}$& $C_{22}$& $C_{23}$ \\
\midrule
Model &&1& -0.05& 0& -1& 0.03& 0 \\\hline
Methods & $\Delta t$ &&&&&&\\\hline
  ODENet &0.001 & 9.93$\times 10^{-1}$ &-5.00$\times 10^{-2}$ &0 &-1.00 &3.00$\times 10^{-2}$& 0\\
 SINDy &0.001 &1.01&-4.99$\times 10^{-2}$&0&-9.95$\times 10^{-1}$  &3.00$\times 10^{-2}$&0\\ \hline
   ODENet&0.01& 9.97$\times 10^{-1}$ &-4.98$\times 10^{-2}$ &0 &-1.00 &2.99$\times 10^{-2}$& 0\\
 SINDy&0.01&1.00&-4.98$\times 10^{-2}$&0&-9.98$\times 10^{-1}$  &3.01$\times 10^{-2}$&0\\ \hline
   ODENet&0.1& 9.82$\times 10^{-1}$ &-4.94$\times 10^{-2}$ &0 &-1.00 &3.01$\times 10^{-2}$& 0\\
 SINDy&0.1&5.87$\times 10^{-1}$&-5.55$\times 10^{-2}$&0&0  &1.83$\times 10^{-1}$&0\\ \hline
   ODENet&0.5& 9.90$\times 10^{-1}$ &-4.99$\times 10^{-2}$ &0 &-1.02 &3.08$\times 10^{-2}$& 0\\
 SINDy&0.5& ---&---&---&---&---&---\\ 
\bottomrule
\bottomrule
\end{tabular}
\caption{Influence of data sampling time step on the accuracy of ODENet and SINDy. $1\%$ white noise is added to the limit cycle case of LV equations. Again, redundant coefficients besides the six ones in the given model \eqref{eq.lv} are omitted.}
\label{tab.lvstep}
\end{table}

Based on the above analysis, the ODENet is applied to learn the dynamics of the LV model \eqref{eq.lv} within different coefficient regimes, see Figure \ref{fig.lv}. As we do not want to introduce any prior knowledge, the right-hand side terms of the ODE system are expanded through a complete polynomial basis. Up to the second order, we have twelve free parameters to learn.  The detailed setup can be found in Box 1.

\begin{figure}[htp]
\begin{tikzpicture}
\node [mybox] (box){%
	\begin{minipage}{0.95\textwidth}
	\textbf{Goal:} Find the correct Lotka-Volterra (LV) model from the given time-series data. \\
	\textbf{Data:} Simulated time trajectories of LV equations representing different types of ODE dynamics in phase space, combined with either small ($1\%$) or large ($10\%$) white noise (with respect to the largest amplitude of data). \\
	\textbf{Setup:} Complete polynomials up to the second order $\bf{\Lambda}=\{1,x_1,x_2,x_1^2,x_1x_2,x_2^2\}$ with twelve adjustable coefficients $\bf{\theta}=(\theta_{ij})_{2\times6}$ are adopted to approximate function $\textbf{f}(\textbf{x})$. For large noise, the noise term $\mathbf{e}(t)$ in \eqref{eq.model} is added as learning parameters too.\\
	\textbf{Learning:} Optimize parameters $\bf{\theta}$ and $\bf{e}$ following the procedure of Algorithm 1. Parameters $\theta_{ij}$ less than the  threshold $\gamma$ is set as zero to remove model redundancy.   Regularization parameter $\mu$ is decreased from $10^{-3}$ to $10^{-5}$, while threshold $\gamma$ is increased from $10^{-4}$ to $10^{-3}$ with iterations. \\
	\textbf{Results:} The correct form of LV models is reproduced with most terms as zeros. The coefficients of remaining terms are close to their expected values, with the maximal relative errors less than $6\%$. The distribution of noise is almost correctly predicted.
	\end{minipage}
};
\node[fancytitle] at (box.north) {\textbf{Box 1: Lotka-Volterra models}};
\end{tikzpicture}
\end{figure}


  
As summarized in Table \ref{tab.lv}, our ODENet shows a competitive performance with the state-of-the-art methods, \textit{e.g.} SINDy \cite{kutz16sparse}.  In fact, for all three typical LV dynamics, all zero terms in the LV model in \eqref{eq.lv} has been correctly picked out by the ODENet through sparse identification. Furthermore, the maximal relative errors between the learned ones and their true values of the remaining nonzero coefficients are less than $6\%$. In contrast, SINDy fails to identify the redundant terms in the spiral and limit-cycle cases, whose performance becomes even worse as the data sampling time steps are increased to $0.01$.

\begin{figure}[!htp]
 \centering
 \includegraphics[width=\textwidth]{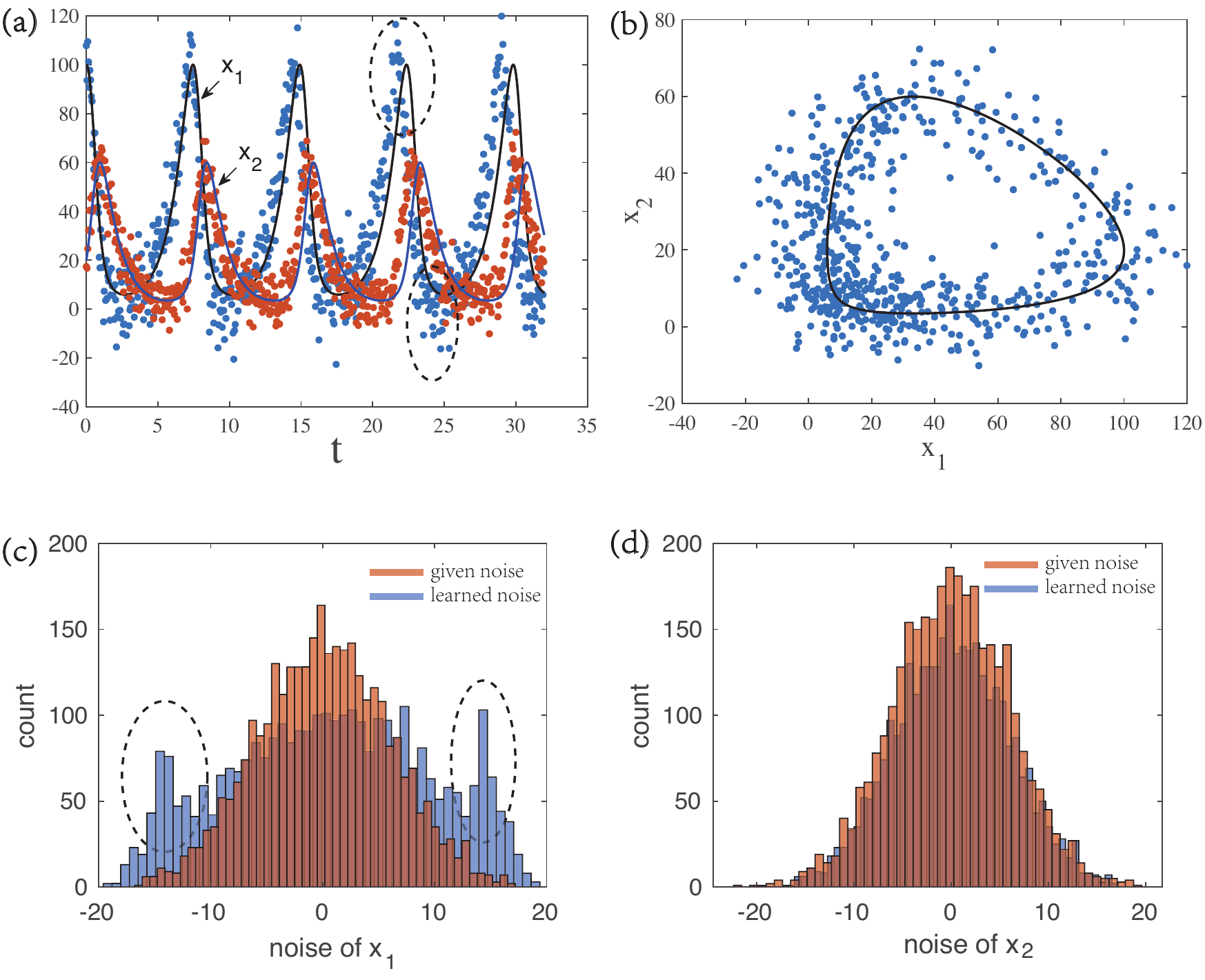}
 \caption{(a-b) Predictions of ODENet on LV equations in the presence of large external noise (up to 10\% of the highest magnitude of the data). Distributions of learned noise in (c) $x_1$ and (d) $x_2$ are compared to the real ones.}
\label{fig.lvnoise}
\end{figure}

Even in the presence of large noise, for example in this case up to $10\%$ white noise with respect to the maximal signal value are added to the data (see Fig. \ref{fig.lvnoise}), our ODENet still shows an astonishing ability in finding out the correct governing equations and revealing the hidden deterministic trajectories which are deeply buried inside noise-spoiled data. The learned noise correctly fits into a Gaussian distribution as expected, though rare events with large displacements are overestimated in the current case, as pointed out in Fig. \ref{fig.lvnoise}c. Further studies show that the deviation from the standard Gaussian distribution disappears as the noise level is lowered (data not shown).  Most importantly, in ODENet, no extra unwanted coefficient will be included in the model as a consequence of sparse identification, even for the flawed and noisy data. This fact is clearly stated through the zero values of $C_{13}$ and $C_{23}$ in the fifth row of Table \ref{tab.lv} for LV equations with large noise.

Compared to difference-based methods, the integration based methods are more tolerant to large time steps in sampling data, as we have claimed. By gradually increasing the sampling time steps of the training data set in the limit-cycle case of LV equations, it is expected that revealing the ODE dynamics from the sample data becomes harder and harder, as less information about the system is included. So that it is not astonishing to see that results of SINDy become untrustable when $\Delta t\geq0.1$. However, our ODENet still works quite well and shows high accuracy in revealing the dynamics even when $\Delta t=0.5$ (see Table \ref{tab.lvstep}).

\subsection{Lorenz equations in chaotic regimes}

\begin{figure}[!htp]
\centering
\includegraphics[width=\linewidth]{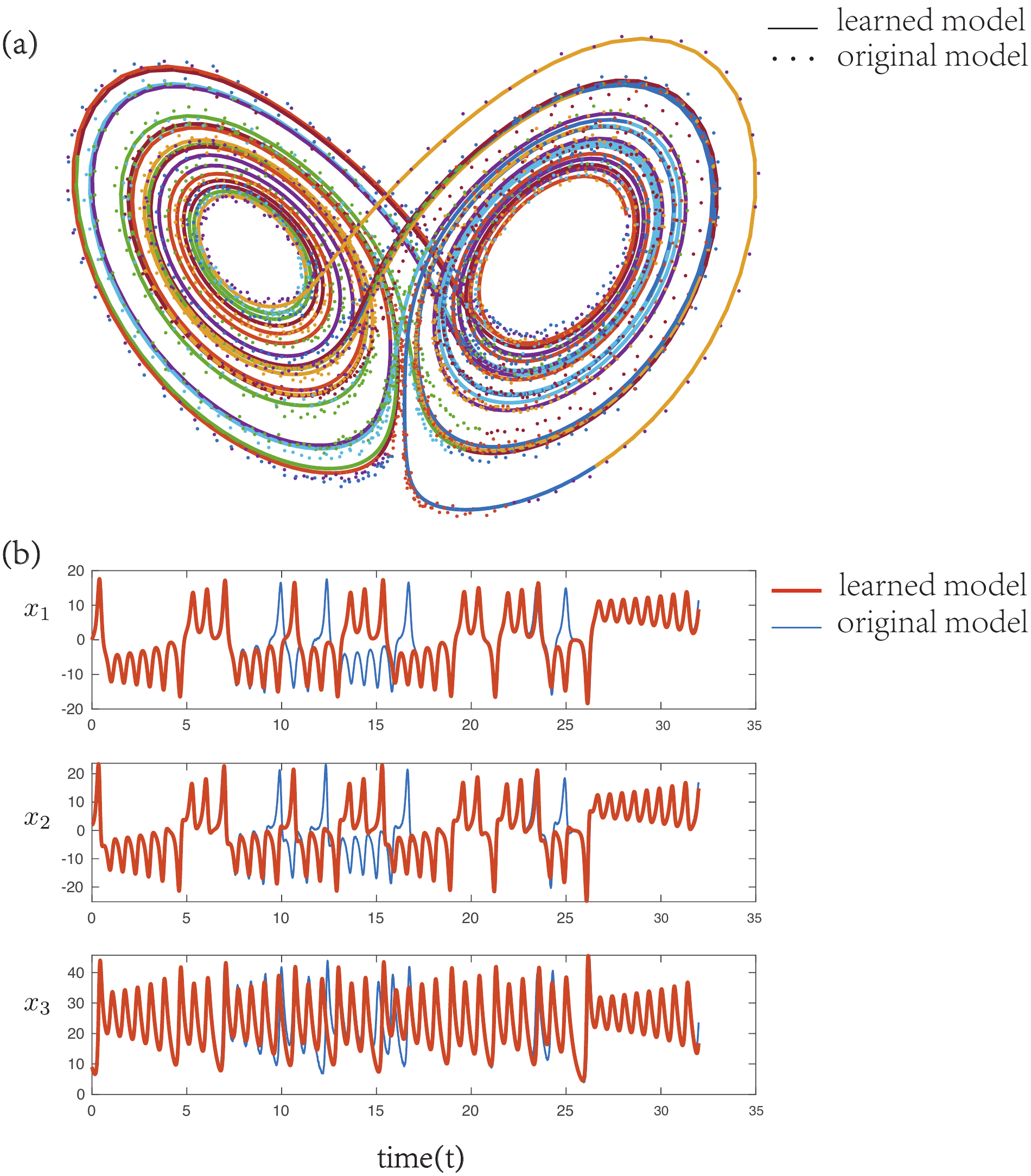}
\caption{Predictions of ODENet on Lorenz equations in the chaotic regime.}
\label{fig.lorenz}
\end{figure}

In the 1960s, American meteorologist Lorenz proposed a simple mathematical system constituted by three ordinary differential equations,
\begin{equation}\label{eq.lorenzp}
\begin{cases}
\frac{\text{d}x_1}{\text{d}t}=C_{11}x_1+C_{12}x_2, \\
\frac{\text{d}x_2}{\text{d}t}=C_{21}x_1+C_{22}x_2+C_{23}x_1x_3, \\
\frac{\text{d}x_3}{\text{d}t}=C_{31}x_3+C_{32}x_1x_2. \\
\end{cases}
\end{equation}
for describing atmospheric turbulence \cite{lorenz1963deterministic}. Lorenz equations became very famous for its chaotic solutions. For a typical parameter combination $C_{11}=-C_{12}=10, C_{21}=28, -C_{32}=8/3, -C_{22}=-C_{23}=C_{31}=1$, the Lorenz system has three equilibrium points, i.e., $\left( {0,0,0} \right)$, $\left( {6\sqrt 2 ,6\sqrt 2 ,27} \right)$ and $\left( { - 6\sqrt 2 , - 6\sqrt 2 ,27} \right)$. Numerical simulation shows that typical trajectories of \eqref{eq.lorenzp} follow a strange attractor in a butterfly shape in the phase space, which first makes a few loops around $\left( {6\sqrt 2 ,6\sqrt 2 ,27} \right)$, then jumps to loops around $\left( { - 6\sqrt 2 , - 6\sqrt 2 ,27} \right)$, and then come back to the point $\left( {6\sqrt 2 ,6\sqrt 2 ,27} \right)$, again and again (see Figure \ref{fig.lorenz}). In this case, solutions of the Lorenz equations are sensitive to  disturbance in the initial conditions, which is widely known as the ``butterfly effect'' in the literature. Therefore, to catch the charming butterfly from chaotic data is an attractive task, which makes the Lorenz system as a benchmark problem for testing the accuracy of numerical schemes as well as the performance of machine learning algorithms.

To increase the learning difficulty, we introduced white noise with magnitude up to $0.5\%$ of the maximal signal data, though no disturbance is included in the initial values. According to the results summarized in Table \ref{tab.lorenz} and Figure \ref{fig.lorenz}, ODENet correctly reproduces the Lorenz attractor, and only a small fraction of trajectories are mispredicted, which is inevitable in the study of chaos.
\\\\
\begin{tikzpicture}
\node [mybox] (box){%
	\begin{minipage}{0.95\textwidth}
	\textbf{Goal:} Predict the strange attractors of Lorenz equations.\\
	\textbf{Data:} Time trajectories of Lorenz equations in the chaotic regime with $0.5\%$ white noise (with respect to the highest amplitude of data) added. \\
	\textbf{Setup:} Complete polynomials up to the second order $\bf{\Lambda}=\{1,x_1,x_2,x_3,x_1^2,x_1x_2,x_1x_3,x_2^2,x_2x_3,x_3^2\}$ with thirty adjustable coefficients $\bf{\theta}=(\theta_{ij})_{3\times10}$ are adopted to approximate function $\textbf{f}(\textbf{x})$. As the orbits are very sensitive to initial values and model coefficients, they are divided into many small pieces to minimize predictive errors.\\
	\textbf{Learning:} Parameters $\bf{\theta}$ are optimized according to the procedure of Algorithm 1. Sparsity requirement is taken.  The regularization factor $\mu$ is set as a decreasing hyper-parameter from $10^{-4}$ to $10^{-8}$ with iterations. The threshold $\gamma$ is an increasing parameter from $10^{-4}$ to $10^{-3}$ to remove redundant terms.
	\\
	\textbf{Results:} The correct form of Lorenz equations are reproduced with the maximal relative errors of coefficients less than $1\%$. The strange attractors are correctly predicted even in a long time.
	\end{minipage}
};
\node[fancytitle] at (box.north) {\textbf{Box 2: Strange attractors of Lorenz equations}};
\end{tikzpicture}

\begin{table}[!htp]\scriptsize
\setlength\extrarowheight{2pt}
\setlength{\tabcolsep}{13.5pt}
\centering
\begin{tabular}{@{} l ccccccc @{}} 
\toprule
\toprule
\multirow{2}*{Lorenz} & \multicolumn{7}{c@{}}{ \hspace{6em} Parameters}\\
\cmidrule(l){2-8}
  & $C_{11}$& $C_{12}$ & $C_{21}$& $C_{22}$& $C_{23}$& $C_{31}$& $C_{32}$ \\
\midrule
Model & -10 & 10 & 28 & -1 & -1 & -8/3 & 1 \\\hline
ODENet& -9.989 & 9.982 & 28.02 & -1.008 & -1.000 & -2.667 & 1.000\\\hline
SINDy & -9.899 & 9.976 & 26.72 & -3.965$\times10^{-1}$ & -9.747$\times10^{-1}$ & -2.447 &9.997$\times10^{-1}$\\
\bottomrule
\bottomrule
\end{tabular}
\caption{Comparison on the accuracy of  ODENet and SINDy for Lorenz equations. For simplicity, only 7 non-zero parameters corresponding to the true model in \eqref{eq.lorenzp} are listed, while the rest 23 parameters are all zeros for ODENet. There is an extra term $-5.115$ in the third equation for SINDy.}
\label{tab.lorenz}
\end{table}

For simplicity, in the above two examples only the second-order complete polynomials have been used to construct the initial model for the ODENet. And luckily the classical LV and Lorenz equations all fall into this category. If higher-order basis functions are adopted, the correct dynamics, in general, could also be revealed, but at a price of larger computational costs and higher risks of facing stiffness problems. This issue has been tested on the LV equations with respect to third-order and fourth-order complete polynomial basis, in which all coefficients for high-order polynomials (above 2) have been correctly eliminated by choosing suitable threshold parameters (data not shown). And in the case of Lorenz equations, the ODENet works in second-order to fifth-order complete polynomial bases but failed to find the correct answer when the right-hand side terms were expanded higher than fifth-order polynomials, probably due to too many redundant terms. We call the attention of readers to the trade-off between the accuracy and efficiency of our method.

\section{Application to the kinetics of actin aggregation}

Actin aggregation into microfilaments is responsible for the contraction of muscle cells and the motility of other cells. In the 1960s, the first analytical molecular model was proposed by Oosawa et al. \cite{oosawa1962theory}, which stated the mechanism of actin aggregation includes three basic steps -- primary nucleation, elongation, and fragmentation. Primary nucleation is an initialization step to generate new growth seeds through a self-organization process. Then small seeds grow into long actin filaments by elongation, meaning monomeric actions are added to the filament ends sequentially. The actin aggregation could be dramatically sped up by fragmentation, through which massive new seeds are generated by breaking long filaments into two shorter pieces without involving primary nucleation. It should be mentioned, besides those forward processes for actin growth, the corresponding inverse processes, like monomer dissociation and fibril annealing, may also make a non-negligible contribution to maintaining the equilibrium distribution of actin filaments. Based on the theory of chemical kinetics, the above picture can be explicitly transformed into a mathematical language of ordinary differential equations, which establishes a direct connection between experimental data and molecular mechanisms of actin growth.

\begin{figure}[htp]
 \centering
 \includegraphics[width=\textwidth]{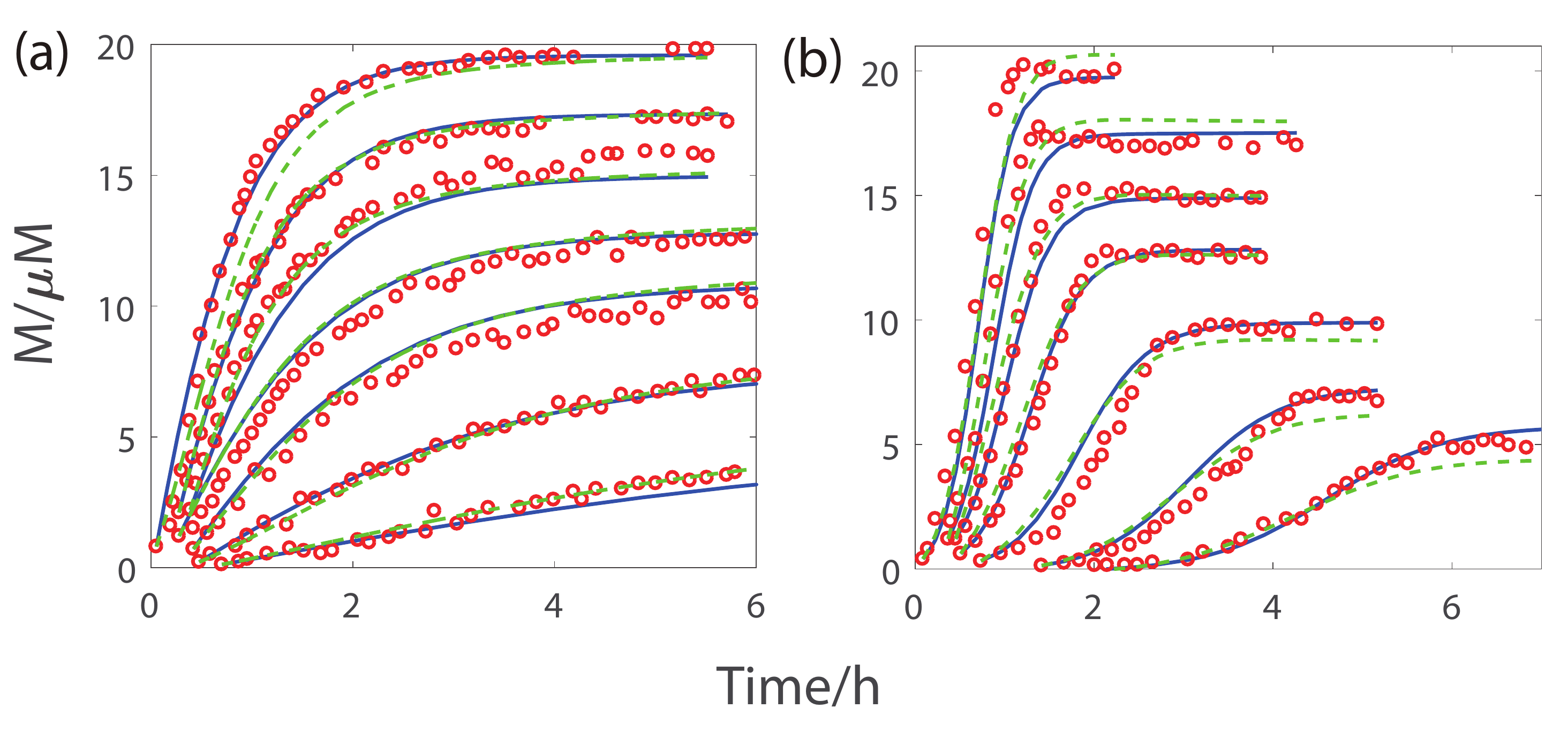}
 \caption{Kinetics of actin aggregation in (a) KCl and (b) MgCl$_2$ solutions respectively. Red circles stand for experimental data in \cite{wegner1982fragmentation}, blue solid lines for predictions of the data-driven model in (\ref{eq.oneequ}), green dashed lines for the physical-based model in (\ref{eq.twoequ}).}
 \label{fig.nu}
\end{figure}

\textit{\textbf {Experimental data: }} In this study, we re-examine the classical experiments done by Wegner et al. \cite{wegner1982fragmentation}, which studied the phenomenon of actin aggregation under two distinct conditions. One is varied concentrations of monomeric actins for $m_{tot}=7.4, 9.6, 12.4, 14.2, 16.2, 18.4, 20.5\mu M$ incubated with $40mM$ \ch{KCl} (Figure \ref{fig.nu}a), the other is $m_{tot}=6.7, 8.5, 11.5, 14.9, 17.3, 20.3, 22.9\mu M$ actins incubated with $0.6mM$ \ch{MgCl_2} and $0.5mM$ EGTA (Figure \ref{fig.nu}b).  The red circles in Figure \ref{fig.nu}a and \ref{fig.nu}b  indicate the mass concentration of actin filaments $M(t)$ which is measured through the ThT fluorescence intensity  under seven different concentrations of monomeric actins.  Clearly, the data are not equally spaced, and the sampling time step $\Delta t$ is not very small. 

To explore the influence of pre-knowledge (or physical insight) on machine learning-based modeling, here we adopt two different setups -- one is purely data-driven, the other is physical-based, which, as we will see, leads to models in distinct forms, but all fit the data quite well.

\textit{\textbf {Purely data-driven model: }}Firstly, we study the pure data-driven modeling without including any pre-knowledge.  To account for the concentration dependence, an additional variable -- the actin monomer concentration $m(t)=m_{tot}-M(t)$  is introduced besides $M(t)$.  As a consequence, we need to learn two ordinary differential equations from the data, i.e.
\begin{equation}\label{eq.oneequ}
\begin{split}
&\frac{\d M}{\d t} = \alpha_0 + \alpha_1 M + \alpha_2m + \alpha_3M^2+ \alpha_4mM+\alpha_5m^2 ,\\
&\frac{\d m}{\d t} = -\alpha_0 - \alpha_1 M - \alpha_2m - \alpha_3M^2- \alpha_4mM-\alpha_5m^2.\\
\end{split}
\end{equation}
Corresponding to reactions up to the second order, terms on the right-hand side are also kept up to the second-order polynomials of $M$ and $m$. Due to the laws of mass conservation, i.e. $M(t)+m(t)=m_{tot}$, five free parameters in the second equation of $m$ can be completely fixed. It is further noted that, since in the current case at least seven concentrations of actin are considered, a global fitting of data with different $m_{tot}$ at the same time is essential for the learning procedure of ODENet.

\begin{table}[]
\setlength{\tabcolsep}{5pt}
 \setlength\extrarowheight{3pt}
\begin{tabular}{ccccccc}
	\toprule
	\toprule
	\multirow{2}{*}{filaments} & \multicolumn{5}{c}{Parameters} & \\ \cline{2-7}
	& ${\alpha _{0}}$ & ${\alpha _{1}}$ & ${\alpha _{2}}$ & ${\alpha _{3}}$ & ${\alpha _{4}}$ & ${\alpha _{5}}$ \\ \hline
	actin in KCl & 4.62$\times10^{-1}$ & -2.16$\times10^{-1}$ & -5.49$\times10^{-1}$ & 5.70$\times10^{-3}$ & 1.10$\times10^{-1}$ & 7.87$\times10^{1}$ \\ \hline
	actin in MgCl$_{2}$ & 0 & -1.41$\times10^{-2}$ & 9.20$\times10^{-3}$ & -3.75$\times10^{-2}$ & 2.28$\times10^{-1}$ & 3.10$\times10^{1}$ \\ 
	\bottomrule
    \bottomrule
\end{tabular}
 \caption{Learned coefficients for data-driven model of actin aggregation.}
\label{tab.re}
\end{table}

Beyond the good agreement between ODENet predictions and experimental data as shown in Figure \ref{fig.nu}, the learned ODE parameters for actin aggregation given in Table \ref{tab.re} are worthy of further clarification, especially their physical meanings. Terms ${\alpha _0}$, ${\alpha _2}m$ and ${\alpha _{\rm{5}}}m^2$ together account for primary nucleation within two monomers. ${\alpha _1}M$ represents the process of degradation (or monomer dissociation). Since it makes a negative contribution to the filament concentration, $\alpha_1$ is always negative as expected. The term ${\alpha_4}mM$ comes from actin filament elongation, which depends on not only the monomer concentration but also the filament concentration. Only the physical meaning of ${\alpha _3}M^2$ is not so straightforward, which may originate from some complicated interactions between filaments, like annealing or clumping. However, based on the coefficients listed in Table \ref{tab.re}, we cannot tell the difference between actin incubating with KCl and with \ch{MgCl_2}. These limitations motivate us to consider a more physical-based model.

\textit{\textbf{Physical based model: }}According to the general theory for actin aggregation \cite{hong2017statistical}, besides the mass concentrations of actin filaments and monomers, the number concentration of actin filaments $P$ also plays a non-negligible role in constructing a complete description of actin growth. So instead of two ODEs, in principle we should consider three coupled equations as the ``correct'' model. However, as the experimental data contain no direct information on $P$, the variable $P$ is actually a hidden one. If we do not write it out explicitly, there is no way to learn it in a purely data-driven modeling. Up to the second-order polynomials, we have
\begin{equation}\label{eq.twoequ}
\begin{split}
&\frac{{dP}}{{dt}} = \alpha_0 +\alpha_1m+\alpha_2m^2 +\alpha_3mM+\alpha_4 P+\alpha_5 M+\alpha_6P^2+\alpha_7PM,\\
&\frac{{dM}}{{dt}} = \alpha_0 +\alpha_1m+2\alpha_2m^2+\alpha_3mM+\alpha_8 M+\alpha_9mP+\alpha_{10}P,\\
&\frac{{dm}}{{dt}} = -\alpha_0 -\alpha_1m-2\alpha_2m^2-\alpha_3mM-\alpha_8 M-\alpha_9mP-\alpha_{10}P.\\
\end{split}
\end{equation}
Here those terms without any physical meaning have been removed, and only eleven free coefficients are kept instead of thirty. The remaining terms all have clear physical interpretations. To be exact, terms $\alpha_0,\alpha_1m$ and $\alpha_2m^2$ stand for primary nucleation, $\alpha_3mM$ for secondary nucleation, $\alpha_4 P$ and $\alpha_8 M$ for fibril degradation, $\alpha_5 M$ for fragmentation, $\alpha_6P^2$ for annealing, $\alpha_7PM$ for clumping, $\alpha_9mP$ for elongation and $\alpha_{10}P$ for monomer dissociation respectively.

 The hidden variable $P$ plays a key role in the physical-based model. But how to learn it is a highly non-trivial task. Since we cannot make a direct comparison between the predicted $P(t)$ with its true values, it will not appear in the loss function. But as $M(t)$ depends on $P(t)$ according to \eqref{eq.twoequ}, optimizing the predictions on $M(t)$ will also lead to an optimization of predicted $P(t)$ simultaneously, as long as the solution of the problem is a fixed-point (it is a general belief, though we are unable to prove its convergence). To help the convergence, the initial values of $P(t)$ were estimated through the approximation $dM/dt\approx\alpha_9mP$, as we regard elongation as the dominant process. Without a proper guess on the initial values, the learning process will fail with a very high probability.

During the learning procedure of ODENet, terms $\alpha_0$, $\alpha_5M$, $\alpha_6P^2$, $\alpha_7PM$, and $\alpha_8M$ are eliminated by sparsity requirement, indicating the corresponding processes may not be essential for modeling. The remaining terms listed in Table \ref{tab.re2} suggest a clear molecular mechanism for the actin aggregation, including primary nucleation (indicated by $\alpha_1m$ and $\alpha_2m^2$), elongation ($\alpha_9mP$), surface catalyzed secondary nucleation ($\alpha_3mM$), monomer dissociation ($\alpha_{10}P$) and degradation $\alpha_4P$. Among them, the first three processes are dominant for microfilament growth, while the latter two are responsible for maintaining the equilibrium state. Further comparing the model coefficients learned from ODENet, the elongation rate for actin aggregation in KCl solution is much smaller than in MgCl$_2$ solution, suggesting the former process is dominated by primary nucleation, while the latter is dominated by elongation and secondary nucleation instead. This dramatic distinction is believed to be caused by different chemical valences of \ch{K^+} and \ch{Mg^{2+}}. Therefore, the physical-based modeling by ODENet indeed provides new insights into those unknown phenomena we are interested in. Qualitatively, the physical-based model uncovered by ODENet is in a good agreement with previous models constructed purely by the human brain \cite{hong2017statistical}. In the latter, primary nucleation, elongation, and fragmentation (a kind of secondary nucleation) were considered as three dominant processes for actin aggregation, and the difference between \ch{KCl} and \ch{MgCl_2} solutions lies in the strength of primary nucleation v.s. secondary nucleation. \\

\begin{table}[]
\setlength{\tabcolsep}{5pt}
\setlength\extrarowheight{5pt}
	\begin{tabular}{ccccccc}
		\toprule
		\toprule
		\multirow{2}{*}{filaments} & \multicolumn{6}{c}{Parameters} \\ \cline{2-7}
		& ${\alpha _{1}}$ & ${\alpha _{2}}$ & ${\alpha _{3}}$ & ${\alpha _{4}}$ & ${\alpha _{9}}$ & ${\alpha _{10}}$ \\ \hline
		actin in KCl & -5.12$\times10^{-2}$ & 7.98$\times10^{-3}$ & 1.16$\times10^{-2}$ & -5.33$\times10^{-1}$ & 7.39$\times10^{-1}$ & -8.82$\times10^{-1}$ \\ \hline
		actin in MgCl$_{2}$ & 2.15$\times10^{-2}$ & 0 & 2.3$\times10^{-2}$ & 0 & 1.19$\times10^{1}$ & -2.97$\times10^{1}$ \\ 
		\bottomrule
		\bottomrule
	\end{tabular}
\caption{Learned coefficients for physical-based model of actin aggregation. Unmentioned coefficients are all zeros.}
\label{tab.re2}
\end{table}

\begin{tikzpicture}
\node [mybox] (box){%
	\begin{minipage}{0.95\textwidth}
	\textbf{Goal:} Compare the purely data-driven model with the physical-based model on kinetics of actin aggregation.\\
	\textbf{Data:} Mass concentration $M(t)$ of actin filaments recorded at different time points in ThT fluorescence experiments. Seven protein concentrations and two buffer conditions are taken into consideration. \\
	\textbf{Setup:} There are two separate setups for the learning procedure:
	\begin{enumerate}
	\item \textbf{Purely data-driven model.} Without any pre-knowledge of the model, we just need single equation of mass concentration $M(t)$ to learn the dynamics. To account for the concentration dependence, an additional variable $m(t)=m_{tot}-M(t)$ is introduced too. Mimicking the function on the right-hand side of ODEs by polynomials up to the second-order, we have to optimize six coefficients $\bf{\theta}=(\theta_{ij})_{2\times 6}$, where $\theta_{2j} = -\theta_{1j},\, j=1,2,\cdots,6$.
	\item \textbf{Physical based model.} According to the general theory for actin aggregation \cite{hong2017statistical}, another hidden variable -- the number concentration of actin filaments $P(t)$ is introduced into the model. Furthermore, we require all kept terms have a clear physical meaning to account for all possible mechanisms for actin growth. In this case, we have three ordinary differential equations with eleven undetermined coefficients. The hidden variable $P(t)$ is constructed in a self-iterative way from the approximation $dM/dt\approx\alpha_9mP$ with a pre-knowledge that elongation makes a major contribution to the mass growth of actin filaments.
	\end{enumerate}
	\textbf{Learning:} Parameters $\bf{\theta}$ are optimized according to the procedure of Algorithm 1. Sparsity requirement is taken.\\
	\textbf{Results:} Two simple models with and without hidden physical variable $P(t)$ are learned separately, both of which can fit ThT trajectories quite well.
	\end{minipage}
};
\node[fancytitle] at (box.north) {\textbf{Box 3: Data-driven v.s. physical-based modeling of actin aggregation}};
\end{tikzpicture}

\section{Conclusion and Discussion}

In this work, we proposed a general and flexible network called ODENet for revealing hidden ODE dynamics from time-series data.  A significant difference of ODENet from the state-of-art regression-based methods like SINDy is the adoption of integration of explicit ODEs along the time trajectory. By further combining with classical machine learning skills, like data batching, back-propagation and optimization,  ODENet inherits the advantages of both machine learning and ODEs. On one hand, the embedding of ODEs makes the whole procedure transparent and interpretable. On the other hand, the schemes of machine learning enable data handling, paralleling, and optimization to be easily and efficiently implemented.  

As illustrated through several novel examples including Lotka-Volterra models for population dynamics, strange attractors of Lorenz equations, and the kinetics of actin aggregation into microfilaments, ODENet shows great merits in several aspects: (1) the ability to deal with data not equally spaced, of a high noise to signal ratio, etc.; (2) tolerance with large sampling time steps; (3) explicitly deriving interpretable models with fewer parameters; (4) efficiently optimizing parameters by BP algorithms; (5) very flexible network structure ready for the incorporation of various new approaches. Therefore, we expect wider applications of ODENet in various branches of natural science, as well as non-trivial extensions to stochastic ODEs and PDEs for a better description of the real world in the near future.

\section*{Acknowledgements}
This work was supported by the National Natural Science Foundation of China (Grant Nos. 21877070 and 11871299) and the Hundred-Talent Program of Sun Yat-Sen University.

\bibliography{ode_modeling}
\bibliographystyle{unsrt}

\end{document}